\documentclass[10pt,notitlepage,reqno]{amsart}%{report}
\usepackage[utf8]{inputenc}

\usepackage{hyperref}
\usepackage[margin=1in]{geometry}
\usepackage{amsmath}
\usepackage{amssymb}
\usepackage[usenames,dvipsnames]{xcolor}
\usepackage[font=footnotesize]{caption}
\usepackage{subcaption}
\usepackage{bbm}
\usepackage{booktabs}
\usepackage{multicol}
\usepackage{arydshln}
\usepackage{enumitem}
\usepackage{setspace}
\usepackage{amsfonts}
\usepackage{mathrsfs}
\usepackage{multirow}
\usepackage[noabbrev,capitalise]{cleveref}
\usepackage{bookmark}
\bookmarksetup{
	numbered,
	open
}
\usepackage{amsthm}
\newtheorem{theorem}{Theorem}[section]
\newtheorem{conjecture}[theorem]{Conjecture}

\theoremstyle{remark}
\newtheorem{remark}[theorem]{Remark}
\numberwithin{equation}{section}
\numberwithin{theorem}{section}
\allowdisplaybreaks
\newcommand{\up}{\big\uparrow}
\newcommand{\down}{\big\downarrow}
\newcommand{\Irr}{\operatorname{Irr}}
\newcommand{\Syl}{\operatorname{Syl}}
\newcommand{\nbhd}{\mathsf{n}}
\newcommand{\G}{\mathbb{G}}
\newcommand{\N}{\mathbb{N}}
\newcommand{\newqed}{\hfill$\lozenge$}

\begin{document}
	
\title{On a conjecture of Hung, Mart\'inez and Navarro}

\author{Stacey Law}
\address{School of Mathematics, Watson Building, University of Birmingham, Edgbaston, Birmingham B15 2TT, UK}
\email{s.law@bham.ac.uk}

\begin{abstract}
	In this note, we prove a recent conjecture of Hung, Mart\'inez and Navarro on the sums of character degrees of finite groups.
\end{abstract}

\keywords{character degree, Sylow normaliser}

\subjclass[2020]{20C15, 20D20}

\maketitle

%========================================================================
\stepcounter{section}

A remarkable development in the representation theory of finite groups was the recent proof of the celebrated McKay Conjecture, completed by Cabanes and Sp\"ath in \cite{CS} following years of contributions by many authors, including a landmark reduction to finite simple groups by Isaacs, Malle and Navarro \cite{IMN}. First posed over fifty years ago, the McKay Conjecture states that for any finite group $G$, any prime $p$ and $P$ a Sylow $p$-subgroup of $G$, there exists a bijection between the set of irreducible characters of $G$ of degree coprime to $p$, and the set of those characters of the Sylow normaliser $N_G(P)$.

A natural question then, is to ask whether such a bijection could be chosen to satisfy additional properties, which would reflect deeper connections between the representation theory of a finite group and its local subgroups. Indeed, many local-global conjectures which refine the McKay Conjecture have been put forward over the years, concerning algebraic properties such as $p$-block structure, Galois actions and character degrees modulo $p$ (see \cite{G} and the references therein); in particular, Giannelli has recently conjectured the following.

\begin{conjecture}[{\cite[Conjecture A]{G}}]\label{conj:G}
	Let $G$ be a finite group, $p$ a prime and $P\in\Syl_p(G)$. Then there exists a bijection
	\[ \Gamma: \Irr_{p'}(G) \longrightarrow \Irr_{p'}(N_G(P)) \]
	such that $\Gamma(\chi)(1)\le\chi(1)$ for all $\chi\in\Irr_{p'}(G)$.
\end{conjecture}

In \cite{HMN}, Hung, Mart\'inez and Navarro reduce Giannelli's strengthening of the McKay Conjecture to a question of finite simple groups, and prove that it holds for all $p=2$. In the course of this work, motivated by Problem 2 from Brauer's famous list of problems \cite{Brauer}, they proposed the following conjecture on the sum of the squares of the degrees of those irreducible characters with $p'$-degree.

\begin{conjecture}[{\cite[Conjecture A]{HMN}}]\label{conj:HMNA}
	Let $G$ be a finite group, $p$ a prime and $P\in\Syl_p(G)$. Then 
	\[ \sum_{\chi\in\Irr_{p'}(G)} \chi(1)^2 \ge |N_G(P):P'| \]
	with equality if and only if $N_G(P)$ has a normal complement in $G$.
\end{conjecture}

They further observe that the inequality part of \cref{conj:HMNA} follows immediately from \cref{conj:G} (so that \cref{conj:HMNA} holds for $p=2$), and then characterise the equality as a consequence of a stronger result concerning normal complements and extendible characters \cite[Theorem 2.6]{HMN}. In this note, we prove that \cref{conj:HMNA} holds for all primes $p$ unconditionally, in particular without needing to assuming \cref{conj:G}.

\begin{proof}[Proof of \cref{conj:HMNA}]
	Let $N:=N_G(P)$.
	First we observe that 
	\[ |N_G(P):P'|=\sum_{\phi\in\Irr_{p'}(N)}\phi(1)^2. \] 
	To see this, note that $|N:P'|=\sum_{\psi\in\Irr(N/P')}\psi(1)^2$. Moreover, from It\^o's Theorem \cite[Theorem 6.15]{IBook}, since $P/P'$ is abelian and normal in $N/P'$ then $\psi(1)\mid|N/P':P/P'|=|N:P|$ for all $\psi\in\Irr(N/P')$. But $p\nmid|N:P|$ so in fact $\Irr(N/P')=\Irr_{p'}(N/P')$. Now we claim that $P'$ lies in the kernel of every $\phi\in\Irr_{p'}(N)$, whence $\Irr_{p'}(N/P')$ can be identified with $\Irr_{p'}(N)$ via lifting (which in particular preserves character degree). But this follows from Clifford's Theorem \cite[Theorem 6.2]{IBook}: given $\phi\in\Irr_{p'}(N)$, then $\phi\down_P=e\sum_{i=1}^t\theta_i$ for some $e\in\N$ where $\{\theta_1,\dotsc,\theta_t\}$ is the $N$-orbit of some $\theta\in\Irr(P)$ and $t=|N:I_N(\theta)|$. But $p\nmid\phi(1)=et\theta(1)$ implies $\theta$ is linear (as $\theta\in\Irr(P)$ gives $\theta(1)\mid|P|$), that is, $\phi\down_P$ is a sum of linear characters of $P$. But $P'$ is contained in the kernel of any linear character of $P$, so $|N_G(P):P'|=\sum_{\phi\in\Irr_{p'}(N)}\phi(1)^2$ as desired.
	
	So it remains to prove that
	\[ \sum_{\chi\in\Irr_{p'}(G)}\chi(1)^2\ge\sum_{\phi\in\Irr_{p'}(N)}\phi(1)^2, \] 
	and to determine the equality condition. Let $X:=\Irr_{p'}(G)$ and $Y:=\Irr_{p'}(N)$. We first define a bipartite graph $\G$ with vertex set $X\sqcup Y$ as follows: the edge set of $\G$ is such that an edge joins $\chi\in X$ and $\phi\in Y$ if and only if $[\chi\down_N,\phi]=[\chi,\phi\up^G]>0$. Given $\chi\in X$, let $\nbhd(\chi) := \{\phi\in\Irr_{p'}(N)\mid [\chi\down_N,\phi]>0 \}$, the neighbourhood of $\chi$ in $\G$. We remark that $\G$ has no isolated vertices: it is clear that every $\chi\in X$ is joined to some $\phi\in Y$ since $p\nmid\chi(1)$. Conversely, given $\phi\in Y$ then $\phi\up^G(1)=|G:N|\phi(1)$, so since $p\nmid|G:N|$ and $p\nmid\phi(1)$ then $\phi\up^G$ has some constituent of $p'$-degree.
	
	Next, we partition $Y$ into subsets as follows. Choose any $\phi_1\in Y$. We have that $\phi_1$ is joined to some $\chi_1\in X$. Set $Y_1:=\nbhd(\chi_1)$.
	If $Y_1\ne Y$, then choose $\phi_2\in Y\setminus Y_1$. Again $\phi_2$ is joined to some $\chi_2\in X$. Moreover, $\chi_2\ne\chi_1$ since by definition $\phi_2\notin Y_1=\nbhd(\chi_1)$. Set $Y_2:=\nbhd(\chi_2)\setminus Y_1$. 
	Continuing, if $Y_1\sqcup Y_2\ne Y$, then choose $\phi_3\in Y\setminus (Y_1\sqcup Y_2)$. There exists $\chi_3\in X$ not equal to $\chi_1,\chi_2$ which is joined to $\phi_3$, and set $Y_3:=\nbhd(\chi_3)\setminus(Y_1\sqcup Y_2)$. 
	Eventually the process terminates as $Y$ is finite, at which point we have found some $m\in\N$ such that we have distinct $\chi_1,\chi_2,\dotsc,\chi_m\in X$ satisfying
	\[ Y=Y_1\sqcup Y_2\sqcup \cdots \sqcup Y_m \quad\text{where}\quad Y_1:=\nbhd(\chi_1),\ Y_i:=\nbhd(\chi_i)\setminus\bigcup_{j<i}Y_j,\ \text{and}\ Y_i\ne\varnothing\quad \forall\ 1\le i\le m. \]
	(In particular, $Y_i\ne\varnothing$ since $\phi_i\in Y_i$ in the notation of our construction, and $m\le|X|$.) Now, clearly for any $i\in\{1,\dotsc,m\}$ we have 
	\[ \chi_i(1) = \sum_{\phi\in\Irr(N)} [\chi_i\down_N,\phi]\phi(1) \ge \sum_{\phi\in\nbhd(\chi_i)}\phi(1) \ge \sum_{\phi\in Y_i}\phi(1), \]
	whence 
	\begin{equation}\label{eq:ineq}
		\sum_{\chi\in\Irr_{p'}(G)}\chi(1)^2 = \sum_{\chi\in X}\chi(1)^2 \ge \sum_{i=1}^m \chi_i(1)^2 \ge \sum_{i=1}^m \sum_{\phi\in Y_i}\phi(1)^2 = \sum_{\phi\in Y}\phi(1)^2 = \sum_{\phi\in\Irr_{p'}(N)}\phi(1)^2,
	\end{equation}
	by convexity of the function $f(x)=x^2$.
	
	Finally, equality in \eqref{eq:ineq} holds if and only if $m=|X|$ and $Y_i=\{\phi_i\}$ and $\chi_i\down_N=\phi_i$ for all $i\in\{1,\dotsc,m\}$.
	Since the (now proven) McKay Conjecture \cite{CS} gives us $|X|=|Y|$,
	this is therefore equivalent to saying that the restriction map $\Irr_{p'}(G)\to\Irr_{p'}(N)$ is a bijection; %without the McKay Conjecture, we'd only know that this map was surjective
	in other words, every $\phi\in\Irr_{p'}(N)$ is extendible to $G$. But this is equivalent to $N_G(P)$ having a normal complement in $G$, by \cite[Theorem 2.6]{HMN} with $N=1$.
\end{proof}

\begin{remark}
	By the same argument, we have that 
	\[ \sum_{\chi\in\Irr_{p'}(G)}\chi(1)^2 \ge \sum_{\phi\in\Irr_{p'}(H)}\phi(1)^2 = |H:P'| \]
	for any subgroup $H\le G$ such that $P\le H\le N_G(P)$. \newqed
\end{remark}

\subsection*{Acknowledgments}
The author thanks Gabriel Navarro for comments on a previous version.


\bigskip

\begin{thebibliography}{9999}
	\bibitem[B63]{Brauer}
	{\sc R.~Brauer,}
	\newblock Representations of finite groups,
	\newblock  pp.~133--175 in \textit{Lectures on modern mathematics}, vol.~1, Wiley, New York, 1963.
	
	\bibitem[CS25]{CS}
	{\sc M.~Cabanes and B.~Sp\"ath,}
	\newblock The McKay Conjecture on character degrees,
	\newblock \textit{Ann.~Math.} \textbf{203} (2026), 933--1032.
	
	\bibitem[G26]{G}
	{\sc E.~Giannelli,}
	\newblock McKay bijections and character degrees,
	\newblock \textit{Adv.~Math.}, \textbf{498} (2026), 111029. %DOI:10.1016/j.aim.2026.111029.
	
	\bibitem[HMN26]{HMN}
	{\sc N.~N.~Hung, J.~M.~Mart\'inez and G.~Navarro,}
	\newblock Sum of the squares of the $p'$-character degrees,
	\newblock \textit{Pac.~J.~Math.} \textbf{342} (2) (2026), 351--380.
	
	\bibitem[I76]{IBook}
	{\sc I.~M.~Isaacs,}
	\newblock \textit{Character theory of finite groups},
	\newblock Dover, New York, 1976.
	
	\bibitem[IMN07]{IMN}
	{\sc I.~M.~Isaacs, G.~Malle and G.~Navarro,}
	\newblock A reduction theorem for the McKay	conjecture,
	\newblock \textit{Invent.~Math.} \textbf{170} (2007), 33–101.
\end{thebibliography}
\end{document}